\documentclass{amsproc}%
\usepackage{graphicx}
\usepackage{amscd}
\usepackage{amsmath}
\usepackage{amsfonts}
\usepackage{amssymb}%
\theoremstyle{plain}

\numberwithin{equation}{section}

\begin{document}
\title[SUPERMATRIX ALGEBRAS]{A NEW\ CLASS\ OF\ SUPERMATRIX ALGEBRAS\ DEFINED\ BY TRANSITIVE MATRICES}
\author{Jen\H{o} Szigeti}
\address{Institute of Mathematics, University of Miskolc, Miskolc, Hungary 3515}
\email{matjeno@uni-miskolc.hu}
\thanks{The author was supported by OTKA K-101515 of Hungary}
\thanks{This research was carried out as part of the TAMOP-4.2.1.B-10/2/KONV-2010-0001
project with support by the European Union, co-financed by the European Social Fund.}

\begin{abstract}
We give a natural definition for the transitivity of a matrix. Using an
endomorphism $\delta:R\longrightarrow R$\ of a base ring $R$ and a transitive
$n\times n$ matrix $T\in\mathrm{M}_{n}(\mathrm{Z}(R))$ over the center
$\mathrm{Z}(R)$, we construct the subalgebra $\mathrm{M}_{n}(R,\delta,T)$\ of
the full $n\times n$ matrix algebra $\mathrm{M}_{n}(R)$ consisting of the so
called $n\times n$ supermatrices. If $\delta^{n}=\mathrm{id}_{R}$\ and $T$
satisfies some extra conditions, then we exhibit an embedding $\overline
{\delta}:R\longrightarrow\mathrm{M}_{n}(R,\delta,T)$. An other result is that
$\mathrm{M}_{n}(R,\delta,T)$\ is closed with respect to taking the
(pre)adjoint. If $R$ is Lie nilpotent and $A\in\mathrm{M}_{n}(R,\delta,T)$,
then the use of the preadjoint and the corresponding determinants and
characteristic polynomials yields a Cayley-Hamilton identity for $A$ with
right coefficients in the fixed ring $\mathrm{Fix}(\delta)$. The presence of a
primitive $n$-th root of unity and $\delta^{n}=\mathrm{id}_{R}$ guarantee the
right integrality of a Lie nilpotent $R$ over $\mathrm{Fix}(\delta)$. We
present essentially new supermatrix algebras over the Grassmann algebra.

\end{abstract}
\maketitle

\noindent1. INTRODUCTION

\bigskip

Throughout the paper a ring $R$\ means a not necessarily commutative ring with
identity, all subrings inherit and all endomorphisms preserve the identity.
The group of units in $R$ is denoted by $\mathrm{U}(R)$ and the centre of $R$
is denoted by $\mathrm{Z}(R)$.

In Section 2 we give a natural definition for the transitivity of an $n\times
n$ matrix over $R$. A complete description of such transitive matrices is
provided. The "blow-up" construction gives an easy way to build bigger
transitive matrices starting from a given one.

In Section 3 we prove that the Hadamard multiplication by a transitive matrix
$T\in\mathrm{M}_{n}(\mathrm{Z}(R))$ gives a new type of automorphisms of the
full $n\times n$ matrix algebra $\mathrm{M}_{n}(R)$. We use an automorphism of
the above type and an endomorphism $\delta_{n}$ of $\mathrm{M}_{n}%
(R)$\ naturally induced by an endomorphism $\delta:R\longrightarrow R$, to
define the subalgebra $\mathrm{M}_{n}(R,\delta,T)$\ of $\mathrm{M}_{n}(R)$
consisting of the so called $n\times n$ supermatrices. If $\delta
^{n}=\mathrm{id}_{R}$\ and $T$ satisfies some extra conditions, then we
exhibit an embedding $\overline{\delta}:R\longrightarrow\mathrm{M}%
_{n}(R,\delta,T)$ of $R$ into the $n\times n$ supermatrix algebra
$\mathrm{M}_{n}(R,\delta,T)$.

Section 4 is devoted to the study of the right and left (and symmetric)
determinants and the corresponding right and left (and symmetric)
characteristic polynomials of supermatrices. The main result of Section 4
claims that the supermatrix algebra $\mathrm{M}_{n}(R,\delta,T)$\ is closed
with respect to taking the preadjoint. As a consequence, we obtain that the
mentioned determinants and the coefficients of the corresponding
characteristic polynomials are in the fixed ring $\mathrm{Fix}(\delta)$ of
$\delta$. If $R$ is Lie nilpotent of index $k$, then we derive that any
supermatrix $A\in\mathrm{M}_{n}(R,\delta,T)$\ satisfies a Cayley-Hamilton
identity (of degree $n^{k}$) with right coefficients in $\mathrm{Fix}(\delta
)$. If $K$ is a field of characteristic zero and there is a primitive $n$-th
root of unity in $K$, then a combination of the embedding result in Section 3
and the Cayley-Hamilton identity in Section 4 gives that a Lie nilpotent
$K$-algebra $R$ is right (and left) integral over $\mathrm{Fix}(\delta)$,
where $\delta:R\longrightarrow R$ is a $K$-automorphism with $\delta
^{n}=\mathrm{id}_{R}$.

In Section 5 we explain in detail how the results of Section 4 generalize the
earlier results in [S2]. In order to demonstrate the importance of
supermatrices, we mention their role in Kemer's theory of T-ideals (see [K]).
New examples of supermatrix algebras over the Grassmann algebra are presented
in 5.2 and 5.3. These examples show that our supermatrix algebras introduced
in Section 3 form a rich class of algebras.

\bigskip

\noindent2. TRANSITIVE\ MATRICES

\bigskip

Let $\mathrm{M}_{n}(R)$\ denote the ring of $n\times n$ matrices over a (not
necessarily commutative) ring $R$ with $1$. A matrix $T=[t_{i,j}]$ in
$\mathrm{M}_{n}(R)$\ is called \textit{transitive} if%
\[
t_{i,i}=1\text{ and }t_{i,j}t_{j,k}=t_{i,k}\text{ for all }i,j,k\in
\{1,\ldots,n\}.
\]
Notice that \textit{ }$t_{i,j}t_{j,i}=t_{i,i}=1$ and $t_{j,i}t_{i,j}%
=t_{j,j}=1$ imply that $t_{i,j}$ and $t_{j,i}$ are (multiplicative) inverses
of each other.

\bigskip

\noindent\textbf{2.1. Proposition.}\textit{ For a matrix }$T\in\mathrm{M}%
_{n}(R)$\textit{\ the following are equivalent.}

\noindent(1)\textit{ }$T$\textit{\ is transitive.}

\noindent(2)\textit{ There exists a sequence }$g_{i}\in\mathrm{U}(R)$\textit{,
}$1\leq i\leq n$\textit{ of invertible elements such that }$t_{i,j}=g_{i}%
g_{j}^{-1}$\textit{ for all }$i,j\in\{1,\ldots,n\}$\textit{. If }$h_{i}%
\in\mathrm{U}(R)$\textit{, }$1\leq i\leq n$\textit{ is an other sequence with
}$t_{i,j}=h_{i}h_{j}^{-1}$\textit{, then }$h_{i}=g_{i}c$\textit{ for some
constant }$c\in\mathrm{U}(R)$\textit{.}

\noindent(3)\textit{ There exists an }$n\times1$\textit{ (column) matrix
}$G=\left[
\begin{array}
[c]{c}%
g_{1}\\
\vdots\\
g_{n}%
\end{array}
\right]  $\textit{ with invertible entries }$g_{i}\in\mathrm{U}(R)$\textit{,
}$1\leq i\leq n$\textit{ such that }$T=G\widetilde{G}$\textit{, where
}$\widetilde{G}=[g_{1}^{-1},\ldots,g_{n}^{-1}]$\textit{ is a }$1\times
n$\textit{ (row) matrix.}

\bigskip

\noindent\textbf{Proof.} (1)$\Longrightarrow$(2): Take $g_{i}=t_{i,1}$, then
the transitivity of $T$ ensures that\textit{ }$t_{i,j}=t_{i,1}t_{1,j}%
=t_{i,1}t_{j,1}^{-1}=g_{i}g_{j}^{-1}$. Clearly, $t_{i,1}=g_{i}g_{1}^{-1}%
=h_{i}h_{1}^{-1}$ implies that $h_{i}=g_{i}c$, where $c=g_{1}^{-1}h_{1}$.

\noindent(2)$\Longrightarrow$(1): Now $t_{i,i}=g_{i}g_{i}^{-1}=1$ and
$t_{i,j}t_{j,k}=g_{i}g_{j}^{-1}g_{j}g_{k}^{-1}=g_{i}g_{k}^{-1}=t_{i,k}$.

\noindent(2)$\Longleftrightarrow$(3): Obvious. $\square$

\bigskip

The \textit{Hadamard product} of the matrices $A=[a_{i,j}]$ and $B=[b_{i,j}]$
in $\mathrm{M}_{n}(R)$ is defined as $A\ast B=[a_{i,j}b_{i,j}]$.

\bigskip

\noindent\textbf{2.2. Proposition.}\textit{ If }$T=[t_{i,j}]$\textit{ is a
transitive matrix in }$\mathrm{M}_{n}(R)$\textit{, then }$T^{2}=nT$\textit{.
If }$R$\textit{ is commutative and }$S\in\mathrm{M}_{n}(R)$\textit{ is an
other transitive matrix, then }$T\ast S$\textit{\ is also transitive.}

\bigskip

\noindent\textbf{Proof.} The $(i,j)$ entry of the square $T^{2}$ is%
\[
\underset{k=1}{\overset{n}{%
{\textstyle\sum}
}}t_{i,k}t_{k,j}=\underset{k=1}{\overset{n}{%
{\textstyle\sum}
}}t_{i,j}=nt_{i,j}.
\]
For a commutative $R$, the transitivity of the Hadamard product $T\ast S$\ is
obvious. $\square$

\bigskip

\noindent\textbf{2.3. Proposition ("blow up").}\textit{ For a transitive
matrix }$T=[t_{i,j}]$\textit{\ in }$\mathrm{M}_{n}(R)$\textit{ and for a
sequence }$0=d_{0}<d_{1}<\cdots<d_{n-1}<d_{n}=m$\textit{ of integers define an
}$m\times m$\textit{ matrix }$\widehat{T}=[\widehat{t}_{p,q}]$\textit{ (the
blow up of }$T$\textit{) as follows:}%
\[
\widehat{t}_{p,q}=t_{i,j}\text{\textit{ if }}d_{i-1}<p\leq d_{i}\text{\textit{
and }}d_{j-1}<q\leq d_{j}\text{.}%
\]
\textit{The above }$\widehat{T}$\textit{ is a transitive matrix in
}$\mathrm{M}_{m}(R)$\textit{. If necessary, we use the notation }%
$T(d_{1},\ldots,d_{n-1},d_{n})$\textit{ instead of }$\widehat{T}$\textit{.}

\bigskip

\noindent\textbf{Proof.} $\widehat{T}=[T_{i,j}]$ can be considered as an
$n\times n$ matrix of blocks, the size of the block $T_{i,j}$\ in the $(i,j)$
position is $(d_{i}-d_{i-1})\times(d_{j}-d_{j-1})$ and each entry of $T_{i,j}$
is $t_{i,j}$. The integers $p,q,r\in\{1,\ldots,m\}$ uniquely determine the
indices $i,j,k\in\{1,\ldots,n\}$ satisfying $d_{i-1}<p\leq d_{i}$,
$d_{j-1}<q\leq d_{j}$ and $d_{k-1}<r\leq d_{k}$. The definition of
$\widehat{T}=[\widehat{t}_{p,q}]$ and the transitivity of $T$ ensure that%
\[
\widehat{t}_{p,q}\widehat{t}_{q,r}=t_{i,j}t_{j,k}=t_{i,k}=\widehat{t}_{p,r}.
\]
Thus the $m\times m$ matrix $\widehat{T}$ is also transitive. $\square$

\bigskip

\noindent\textbf{2.4. Examples.} For an invertible element $u\in\mathrm{U}%
(R)$, the sequence $g_{i}=u^{i-1}$, $1\leq i\leq n$ in Proposition 2.1 gives
an $n\times n$ matrix $P^{(u)}=[p_{i,j}]$ with $p_{i,j}=u^{i-j}$. The choice
$u=1$ yields the Hadamard identity $H_{n}$ (each entry of $H_{n}$ is $1$). If
$n=2$ and $u=-1$, then we obtain the following $2\times2$ matrix%
\[
P^{(-1)}=P=\left[
\begin{array}
[c]{cc}%
1 & -1\\
-1 & 1
\end{array}
\right]  .
\]
For $d_{1}=d$ and $d_{2}=m$, the blow up%
\[
\widehat{P}=P(d,m)=\left[
\begin{array}
[c]{cc}%
P_{1,1} & P_{1,2}\\
P_{2,1} & P_{2,2}%
\end{array}
\right]  ,
\]
of $P$ (see Proposition 2.3) contains the square blocks $P_{1,1}$ and
$P_{2,2}$ of sizes $d\times d$ and $(m-d)\times(m-d)$ and the rectangular
blocks $P_{1,2}$ and $P_{2,1}$ of sizes $d\times(m-d)$ and $(m-d)\times d$.
Each entry of $P_{1,1}$ and $P_{2,2}$ is $1$ and each entry of $P_{1,2}$ and
$P_{2,1}$ is $-1$. Thus $P^{(u)}$, $H_{n}$, $P$ and $P(d,m)$ are examples of
transitive matrices. $\square$

\bigskip

\noindent3. THE\ ALGEBRA\ OF\ SUPERMATRICES

\bigskip

In the present section we consider certain endomorphisms of the full matrix
algebra $\mathrm{M}_{n}(R)$. A typical example is the conjugate automorphism
$X\longmapsto W^{-1}XW$, where $W\in\mathrm{GL}_{n}(R)$ is an invertible
matrix (see the well known Skolem-Noether theorem). Any endomorphism
$\delta:R\longrightarrow R$ of $R$ can be naturally extended to an
endomorphism $\delta_{n}:\mathrm{M}_{n}(R)\longrightarrow\mathrm{M}_{n}(R)$.
The following proposition provides a further type of automorphisms of
$\mathrm{M}_{n}(R)$.

\bigskip

\noindent\textbf{3.1. Proposition.}\textit{ Let }$T=[t_{i,j}]$\textit{ be a
matrix in }$\mathrm{M}_{n}(\mathrm{Z}(R))$\textit{, where }$\mathrm{Z}%
(R)$\textit{ denotes the center of }$R$\textit{. The following conditions are
equivalent:}

\noindent(1)\textit{ }$T$\textit{ is transitive,}

\noindent(2)\textit{ for }$A\in\mathrm{M}_{n}(R)$\textit{\ the map (Hadamard
multiplication by }$T$\textit{) }$\Theta_{T}(A)=T\ast A$\textit{ is an
automorphism of the matrix algebra }$\mathrm{M}_{n}(R)$\textit{ (over
}$\mathrm{Z}(R)$\textit{).}

\bigskip

\noindent\textbf{Proof.} (1)$\Longrightarrow$(2): In order to prove the
multiplicative property of $\Theta_{T}$, it is enough to check that
$T\ast(AB)=(T\ast A)(T\ast B)$ for all $A,B\in\mathrm{M}_{n}(R)$. Indeed, the
$(i,j)$ entries of $(T\ast A)(T\ast B)$ and $T\ast(AB)$ are equal:%
\[
\underset{k=1}{\overset{n}{%
{\textstyle\sum}
}}t_{i,k}a_{i,k}t_{k,j}b_{k,j}=\underset{k=1}{\overset{n}{%
{\textstyle\sum}
}}t_{i,k}t_{k,j}a_{i,k}b_{k,j}=\underset{k=1}{\overset{n}{%
{\textstyle\sum}
}}t_{i,j}a_{i,k}b_{k,j}=t_{i,j}\underset{k=1}{\overset{n}{%
{\textstyle\sum}
}}a_{i,k}b_{k,j}.
\]
The inverse of $\Theta_{T}$ is $\Theta_{T}^{-1}(A)=S\ast A$, where
$S=[t_{i,j}^{-1}]$ is also transitive in $\mathrm{M}_{n}(\mathrm{Z}(R))$.

\noindent(2)$\Longrightarrow$(1): Now $t_{i,i}=1$ is a consequence of $T\ast
I_{n}=I_{n}$. Using the standard matrix units $E_{i,j}$ and $E_{j,k}$ in
$\mathrm{M}_{n}(R)$, the multiplicative property of $\Theta_{T}$\ gives that%
\[
t_{i,k}E_{i,k}=T\ast E_{i,k}=T\ast(E_{i,j}E_{j,k})=(T\ast E_{i,j})(T\ast
E_{j,k})
\]%
\[
=(t_{i,j}E_{i,j})(t_{j,k}E_{j,k})=t_{i,j}t_{j,k}E_{i,k},
\]
whence $t_{i,k}=t_{i,j}t_{j,k}$ follows. $\square$

\bigskip

If $\Delta,\Theta:S\longrightarrow S$\ are endomorphisms of the ring $S$, then
the subset%
\[
S(\Delta=\Theta)=\{s\in S\mid\Delta(s)=\Theta(s)\}
\]
is a subring of $S$ (notice that $\Delta(1)=\Theta(1)=1$). Let $\mathrm{Fix}%
(\Delta)=S(\Delta=\mathrm{id}_{S})$ denote the subring of the fixed elements
of $\Delta$.

Now take $S=\mathrm{M}_{n}(R)$, $\Delta=\delta_{n}$ and $\Theta_{T}(A)=T\ast
A$, where $\delta:R\longrightarrow R$ is an endomorphism and $T=[t_{i,j}]$ is
a transitive matrix in $\mathrm{M}_{n}(\mathrm{Z}(R))$. The short notation for
$\mathrm{M}_{n}(R)(\delta_{n}=\Theta_{T})$ is%
\[
\mathrm{M}_{n}(R,\delta,T)=\{A\in\mathrm{M}_{n}(R)\mid A=[a_{i,j}]\text{ and
}\delta(a_{i,j})=t_{i,j}a_{i,j}\text{ for all }1\leq i,j\leq n\}.
\]
If $C\subseteq\mathrm{Z}(R)\cap\mathrm{Fix}(\delta)$ is a (commutative)
subring (say $C=\mathbb{Z}$), then $\mathrm{M}_{n}(R,\delta,T)$ is a
$C$-subalgebra of $\mathrm{M}_{n}(R)$. The elements of the \textit{supermatrix
algebra} $\mathrm{M}_{n}(R,\delta,T)$\ are called $(\delta,T)$%
\textit{-supermatrices}. If $t_{i,j}\in\mathrm{Fix}(\delta)$ for all $1\leq
i,j\leq n$, then $\mathrm{M}_{n}(R,\delta,T)$ is closed with respect to the
action of $\delta_{n}$.

\bigskip

\noindent\textbf{3.2. Theorem.}\textit{ Let }$T=[t_{i,j}]$\textit{ be an
}$n\times n$\textit{\ transitive matrix such that the entries }$t_{i,1}%
\in\mathrm{U}(R)\cap\mathrm{Z}(R)$\textit{, }$1\leq i\leq n$\textit{ of the
first column are central invertible elements. If }$\frac{1}{n}\in R$\textit{
and }$\delta:R\longrightarrow R$\textit{ is an arbitrary endomorphism, then
for }$r\in R$\textit{\ take }$\overline{\delta}(r)=\frac{1}{n}\left[
x_{i,j}(r)\right]  _{n\times n}$\textit{, where}%
\[
x_{i,j}(r)=\underset{k=0}{\overset{n-1}{%
{\textstyle\sum}
}}t_{j,i}^{k}\delta^{k}(r)=\underset{k=0}{\overset{n-1}{%
{\textstyle\sum}
}}t_{i,j}^{-k}\delta^{k}(r)
\]
\textit{and }$\frac{1}{n}x_{i,j}(r)$\textit{\ is in the }$(i,j)$\textit{
position of the }$n\times n$\textit{ matrix }$\overline{\delta}(r)$\textit{.
The above definition gives a map }$\overline{\delta}:R\longrightarrow
\mathrm{M}_{n}(R)$\textit{ such that }$\overline{\delta}(cr)=c\overline
{\delta}(r)$\textit{ and }$\overline{\delta}(rc)=\overline{\delta}%
(r)c$\textit{ for all }$c\in\mathrm{Fix}(\delta)$\textit{.}

\noindent(1)\textit{ If }$t_{i,1}^{n}=1$\textit{ and }$1-t_{i,j}$\textit{ is a
non-zero divisor in }$\mathrm{Z}(R)$\textit{ for all }$1\leq i,j\leq
n$\textit{ with }$i\neq j$\textit{, then }$\overline{\delta}(r)=rI_{n}%
$\textit{ is a scalar matrix for any }$r\in\mathrm{Fix}(\delta)$\textit{.}

\noindent(2)\textit{ If }$t_{1,1}^{k}+t_{2,1}^{k}+\cdots+t_{n,1}^{k}%
=t_{1,1}^{-k}+t_{2,1}^{-k}+\cdots+t_{n,1}^{-k}=0$\textit{\ for all }$1\leq
k\leq n-1$\textit{, then }$\overline{\delta}:R\longrightarrow\mathrm{M}%
_{n}(R)$\textit{ is an embedding of rings.}

\noindent(3)\textit{ If }$t_{i,1}\in\mathrm{Fix}(\delta)$\textit{, }%
$t_{i,1}^{n}=1$\textit{ for all }$1\leq i\leq n$ \textit{and }$\delta
^{n}=\mathrm{id}_{R}$\textit{, then }$\overline{\delta}$\textit{ is an
}$R\longrightarrow\mathrm{M}_{n}(R,\delta,T)$\textit{\ embedding.}

\bigskip

\noindent\textbf{Proof.} Clearly, $t_{i,j}=t_{i,1}t_{j,1}^{-1}$ ensures that
$t_{i,j}\in\mathrm{U}(R)\cap\mathrm{Z}(R)$ for all $1\leq i,j\leq n$. Since
$\delta^{k}(cr)=c\delta^{k}(r)$ and $\delta^{k}(rc)=\delta^{k}(r)c$ hold for
all $r\in R$ and $c\in\mathrm{Fix}(\delta)$, we deduce that $\overline{\delta
}(cr)=c\overline{\delta}(r)$ and $\overline{\delta}(rc)=\overline{\delta}(r)c$.

\noindent(1) In view of $t_{i,i}=1$, the $(i,i)$\ diagonal entry of
$\overline{\delta}(r)$\ is%
\[
\overline{\delta}(r)_{i,i}=\frac{1}{n}\underset{k=0}{\overset{n-1}{%
{\textstyle\sum}
}}t_{i,i}^{k}\delta^{k}(r)=\frac{1}{n}\underset{k=0}{\overset{n-1}{%
{\textstyle\sum}
}}r=\frac{1}{n}(nr)=r.
\]
For $i\neq j$, the $(i,j)$\ non-diagonal entry of $\overline{\delta}(r)$\ is%
\[
\overline{\delta}(r)_{i,j}=\frac{1}{n}\underset{k=0}{\overset{n-1}{%
{\textstyle\sum}
}}t_{j,i}^{k}\delta^{k}(r)=\frac{1}{n}\left(  \underset{k=0}{\overset{n-1}{%
{\textstyle\sum}
}}t_{j,i}^{k}\right)  r
\]
and%
\[
\underset{k=0}{\overset{n-1}{%
{\textstyle\sum}
}}t_{j,i}^{k}=1+t_{j,i}+t_{j,i}^{2}+\cdots+t_{j,i}^{n-1}=0
\]
follows from $t_{j,i}^{n}=t_{j,1}^{n}t_{i,1}^{-n}=1$ and%
\[
(1-t_{j,i})\left(  1+t_{j,i}+t_{j,i}^{2}+\cdots+t_{j,i}^{n-1}\right)
=1-t_{j,i}^{n}=0.
\]
\noindent(2) The additive property of $\overline{\delta}$ is clear. In order
to prove the multiplicative property of $\overline{\delta}$, take $r,s\in R$
and compute the $(i,j)$ entry in the product of the $n\times n$ matrices
$\overline{\delta}(r)$ and $\overline{\delta}(s)$:%
\[
\left(  \overline{\delta}(r)\overline{\delta}(s)\right)  _{i,j}=\frac{1}%
{n^{2}}\underset{u=1}{\overset{n}{\sum}}x_{i,u}(r)x_{u,j}(s)=\frac{1}{n^{2}%
}\underset{u=1}{\overset{n}{\sum}}\left(  \underset{k=0}{\overset{n-1}{%
{\textstyle\sum}
}}t_{u,i}^{k}\delta^{k}(r)\right)  \left(  \underset{l=0}{\overset{n-1}{%
{\textstyle\sum}
}}t_{j,u}^{l}\delta^{l}(s)\right)
\]%
\[
=\frac{1}{n^{2}}\underset{u=1}{\overset{n}{\sum}}\left(  \underset{0\leq
k,l\leq n-1}{\sum}t_{u,i}^{k}\delta^{k}(r)t_{j,u}^{l}\delta^{l}(s)\right)
=\frac{1}{n^{2}}\underset{0\leq k,l\leq n-1}{\sum}\left(  \underset
{u=1}{\overset{n}{\sum}}t_{j,u}^{l}t_{u,i}^{k}\right)  \delta^{k}(r)\delta
^{l}(s)
\]%
\[
\overset{(\ast)}{=}\frac{1}{n^{2}}\underset{0\leq q\leq n-1}{\sum}\left(
\underset{u=1}{\overset{n}{\sum}}t_{j,u}^{q}t_{u,i}^{q}\delta^{q}(r)\delta
^{q}(s)\right)  =\frac{1}{n^{2}}\underset{0\leq q\leq n-1}{\sum}nt_{j,i}%
^{q}\delta^{q}(r)\delta^{q}(s)
\]%
\[
=\frac{1}{n}\underset{q=0}{\overset{n-1}{\sum}}t_{j,i}^{q}\delta^{q}%
(rs)=\frac{1}{n}x_{i,j}(rs)=\left(  \overline{\delta}(rs)\right)  _{i,j}.
\]
The essential part of the above calculation is step $(\ast)$. We used that%
\[
\underset{u=1}{\overset{n}{\sum}}t_{j,u}^{l}t_{u,i}^{k}=\underset
{u=1}{\overset{n}{\sum}}(t_{j,1}t_{u,1}^{-1})^{l}(t_{u,1}t_{1,i})^{k}%
=t_{j,1}^{l}\left(  \underset{u=1}{\overset{n}{\sum}}t_{u,1}^{k-l}\right)
t_{1,i}^{k}=0
\]
for any fixed pair $(k,l)$ with $0\leq k,l\leq n-1$ and $k\neq l$. If
$\overline{\delta}(r)=0_{n\times n}$, then%
\[
0=\underset{i=1}{\overset{n}{%
{\textstyle\sum}
}}x_{1,i}(r)=\underset{i=1}{\overset{n}{%
{\textstyle\sum}
}}\left(  \underset{k=0}{\overset{n-1}{%
{\textstyle\sum}
}}t_{i,1}^{k}\delta^{k}(r)\right)  =nr+\underset{k=1}{\overset{n-1}{%
{\textstyle\sum}
}}\left(  \underset{i=1}{\overset{n}{%
{\textstyle\sum}
}}t_{i,1}^{k}\right)  \delta^{k}(r)=nr,
\]
whence $r=\frac{1}{n}(nr)=0$ follows. Thus $\overline{\delta}$ is an embedding
of rings.

\noindent(3) Now $t_{j,i}^{n}=t_{j,1}^{n}t_{i,1}^{-n}=1$ and $\delta
^{n}=\mathrm{id}_{R}$ imply that $t_{j,i}^{n}\delta^{n}(r)=r=t_{j,i}^{0}%
\delta^{0}(r)$. Since $\delta(t_{j,i}^{k}\delta^{k}(r))=t_{j,i}^{k}%
\delta^{k+1}(r)$ is a consequence of $t_{j,i}^{k}=t_{j,1}^{k}t_{i,1}^{-k}%
\in\mathrm{Fix}(\delta)$, we have%
\[
\delta(x_{i,j}(r))=\underset{k=0}{\overset{n-1}{%
{\textstyle\sum}
}}t_{j,i}^{k}\delta^{k+1}(r)=\underset{l=1}{\overset{n}{%
{\textstyle\sum}
}}t_{j,i}^{l-1}\delta^{l}(r)=t_{i,j}\underset{l=1}{\overset{n}{%
{\textstyle\sum}
}}t_{j,i}^{l}\delta^{l}(r)=t_{i,j}x_{i,j}(r)
\]
proving that $\overline{\delta}(r)\in\mathrm{M}_{n}(R,\delta,T)$. $\square$

\bigskip

\noindent\textbf{3.3. Remark.} In the presence of $t_{1,1}^{n}=t_{2,1}%
^{n}=\cdots=t_{n,1}^{n}$ condition $t_{1,1}^{k}+t_{2,1}^{k}+\cdots+t_{n,1}%
^{k}=0$ implies $t_{1,1}^{k-n}+t_{2,1}^{k-n}+\cdots+t_{n,1}^{k-n}=0$. Thus
$t_{1,1}^{-k}+t_{2,1}^{-k}+\cdots+t_{n,1}^{-k}=0$ in Theorem 3.2 is
superfluous if $t_{1,1}^{n}=t_{2,1}^{n}=\cdots=t_{n,1}^{n}$.

\bigskip

\noindent\textbf{3.4. Proposition.}\textit{ Let }$R$\textit{ be an algebra
over a field }$K$\textit{ of characteristic zero (notice that }$K\subseteq
\mathrm{Z}(R)$\textit{). If }$e\in K$\textit{ is a primitive }$n$%
\textit{-th\ root of unity (}$e^{n}=1\neq e^{k}$\textit{ for all }$1\leq k\leq
n-1$\textit{), then for the }$n\times n$\textit{\ transitive matrix }%
$P^{(e)}=[p_{i,j}]$\textit{ with }$p_{i,j}=e^{i-j}$\textit{, }$1\leq i,j\leq
n$\textit{ (see Example }2.4\textit{) we have }$p_{i,1}^{n}=1$\textit{,
}$p_{1,1}^{k}+p_{2,1}^{k}+\cdots+p_{n,1}^{k}=0$\textit{ and }$1-p_{i,j}%
$\textit{ is a non-zero divisor (invertible) in }$K$\textit{ for all }$1\leq
i,j\leq n$\textit{ with }$i\neq j$\textit{.}

\bigskip

\noindent\textbf{Proof.} Now $p_{i,1}^{n}=(e^{i-1})^{n}=(e^{n})^{i-1}=1$. If
$1\leq k\leq n-1$, then $1-e^{k}\neq0$ is invertible in $K$ and%
\[
(1-e^{k})(1+e^{k}+e^{2k}+\cdots+e^{(n-1)k})=1-e^{nk}=0
\]
implies that%
\[
p_{1,1}^{k}+p_{2,1}^{k}+\cdots+p_{n,1}^{k}=1+e^{k}+e^{2k}+\cdots
+e^{(n-1)k}=0.
\]
If $i\neq j$, then $1-p_{i,j}=1-e^{i-j}\neq0$ is invertible in $K$. $\square$

\bigskip

\noindent\textbf{3.5. Corollary.}\textit{ Let }$\frac{1}{2}\in R$\textit{ and
}$\delta:R\longrightarrow R$\textit{ be an arbitrary endomorphism. For }$r\in
R$\textit{ the definition}%
\[
\overline{\delta}(r)=\frac{1}{2}\left[
\begin{array}
[c]{cc}%
r+\delta(r) & r-\delta(r)\\
r-\delta(r) & r+\delta(r)
\end{array}
\right]
\]
\textit{gives an embedding }$\overline{\delta}:R\longrightarrow\mathrm{M}%
_{2}(R)$\textit{. If }$\delta^{2}=\mathrm{id}_{R}$\textit{, then }%
$\overline{\delta}$\textit{ is an }$R\longrightarrow\mathrm{M}_{2}%
(R,\delta,P)$\textit{\ embedding (for }$P$\textit{\ see }2.4\textit{).}

\bigskip

\noindent\textbf{Proof.} Take $n=2$ and $t_{1,1}=t_{2,2}=1$, $t_{1,2}%
=t_{2,1}=-1$ in Theorem 3.2. $\square$

\bigskip

\noindent4. THE\ RIGHT\ AND\ LEFT\ DETERMINANTS OF\ A SUPERMATRIX

\bigskip

The following definitions and the basic results about the symmetric and
Lie-nilpotent analogues of the classical determinant theory can be found in
[Do, S1, S3, SvW].

Let $\mathrm{S}_{n}$ denote the symmetric group of all permutations of the set
$\{1,2,\ldots,n\}$. For an $n\times n$ matrix $A=[a_{i,j}]$ over an arbitrary
(possibly non-commutative) ring or algebra $R$ with $1$, the element%
\[
\mathrm{sdet}(A)=\underset{\tau,\pi\in\mathrm{S}_{n}}{\sum}\mathrm{sgn}%
(\pi)a_{\tau(1),\pi(\tau(1))}\cdots a_{\tau(t),\pi(\tau(t))}\cdots
a_{\tau(n),\pi(\tau(n))}%
\]%
\[
=\underset{\alpha,\beta\in\mathrm{S}_{n}}{\sum}\mathrm{sgn}(\alpha
)\mathrm{sgn}(\beta)a_{\alpha(1),\beta(1)}\cdots a_{\alpha(t),\beta(t)}\cdots
a_{\alpha(n),\beta(n)}%
\]
of $R$ is the \textit{symmetric determinant} of $A$. The \textit{preadjoint
matrix} $A^{\ast}=[a_{r,s}^{\ast}]$ of $A=[a_{i,j}]$ is defined as the
following natural symmetrization of the classical adjoint:%
\[
a_{r,s}^{\ast}=\underset{\tau,\pi}{\sum}\mathrm{sgn}(\pi)a_{\tau(1),\pi
(\tau(1))}\cdots a_{\tau(s-1),\pi(\tau(s-1))}a_{\tau(s+1),\pi(\tau
(s+1))}\cdots a_{\tau(n),\pi(\tau(n))}%
\]%
\[
=\underset{\alpha,\beta}{\sum}\mathrm{sgn}(\alpha)\mathrm{sgn}(\beta
)a_{\alpha(1),\beta(1)}\cdots a_{\alpha(s-1),\beta(s-1)}a_{\alpha
(s+1),\beta(s+1)}\cdots a_{\alpha(n),\beta(n)}\text{ },
\]
where the first sum is taken over all $\tau,\pi\in\mathrm{S}_{n}$ with
$\tau(s)=s$ and $\pi(s)=r$ (while the second sum is taken over all
$\alpha,\beta\in\mathrm{S}_{n}$ with $\alpha(s)=s$ and $\beta(s)=r$). We note
that the $(r,s)$ entry of $A^{\ast}$ is exactly the signed symmetric
determinant $(-1)^{r+s}\mathrm{sdet}(A_{s,r})$\ of the $(n-1)\times
(n-1)$\ minor $A_{s,r}$\ of $A$ arising from the deletion of the $s$-th row
and the $r$-th column of $A$. If $R$\ is commutative, then $\mathrm{sdet}%
(A)=n!\mathrm{\det}(A)$ and $A^{\ast}=(n-1)!\mathrm{adj}(A)$, where
$\mathrm{\det}(A)$ and $\mathrm{adj}(A)$ denote the ordinary determinant and
adjoint of $A$.

The right adjoint sequence $(P_{k})_{k\geq1}$ of $A$ is defined by the
recursion: $P_{1}=A^{\ast}$ and $P_{k+1}=(AP_{1}\cdots P_{k})^{\ast}$ for
$k\geq1$. The $k$-th right determinant is the trace of $AP_{1}\cdots P_{k}$:%
\[
\mathrm{rdet}_{(k)}(A)=\mathrm{tr}(AP_{1}\cdots P_{k}).
\]
The left adjoint sequence $(Q_{k})_{k\geq1}$ can be defined analogously:
$Q_{1}=A^{\ast}$ and $Q_{k+1}=(Q_{k}\cdots Q_{1}A)^{\ast}$ for $k\geq1$. The
$k$-th left determinant of $A$ is%
\[
\mathrm{ldet}_{(k)}(A)=\mathrm{tr}(Q_{k}\cdots Q_{1}A).
\]
Clearly, $\mathrm{rdet}_{(k+1)}(A)=\mathrm{rdet}_{(k)}(AA^{\ast})$ and
$\mathrm{ldet}_{(k+1)}(A)=\mathrm{ldet}_{(k)}(A^{\ast}A)$. We note that%
\[
\mathrm{rdet}_{(1)}(A)=\mathrm{tr}(AA^{\ast})=\mathrm{sdet}(A)=\mathrm{tr}%
(A^{\ast}A)=\mathrm{ldet}_{(1)}(A).
\]
As we can see in Section 5, the following theorem is a broad generalization of
one of the main results in [S2].

\bigskip

\noindent\textbf{4.1. Theorem.}\textit{ Let }$\delta:R\longrightarrow
R$\textit{ be an endomorphism and }$T=[t_{i,j}]$\textit{ be a transitive
matrix in }$\mathrm{M}_{n}(\mathrm{Z}(R))$\textit{. If }$A\in\mathrm{M}%
_{n}(R,\delta,T)$\textit{ is a supermatrix, then }$A^{\ast}\in\mathrm{M}%
_{n}(R,\delta,T)$\textit{. In other words, the supermatrix algebra
}$\mathrm{M}_{n}(R,\delta,T)$\textit{\ is closed with respect to taking the
preadjoint.}

\bigskip

\noindent\textbf{Proof.} The $(r,s)$ entry of $A^{\ast}$ is%
\[
a_{r,s}^{\ast}=\underset{\tau,\pi}{\sum}\mathrm{sgn}(\pi)a_{\tau(1),\pi
(\tau(1))}\cdots a_{\tau(s-1),\pi(\tau(s-1))}a_{\tau(s+1),\pi(\tau
(s+1))}\cdots a_{\tau(n),\pi(\tau(n))},
\]
where $A=[a_{i,j}]$ and the sum is taken over all $\tau,\pi\in\mathrm{S}_{n}$
with $\tau(s)=s$ and $\pi(s)=r$. In order to see that $A^{\ast}\in
\mathrm{M}_{n}(R,\delta,T)$, we prove that $\delta(a_{r,s}^{\ast}%
)=t_{r,s}a_{r,s}^{\ast}$ for all $1\leq r,s\leq n$. Using $\delta
(a_{\tau(i),\pi(\tau(i))})=t_{\tau(i),\pi(\tau(i))}a_{\tau(i),\pi(\tau(i))}$
we obtain that%
\[
\delta(a_{r,s}^{\ast})=\underset{\tau,\pi}{\sum}\mathrm{sgn}(\pi)b(1,\tau
,\pi)\cdots b(s-1,\tau,\pi)b(s+1,\tau,\pi)\cdots b(n,\tau,\pi),
\]
where $b(i,\tau,\pi)=t_{\tau(i),\pi(\tau(i))}a_{\tau(i),\pi(\tau(i))}$. Since
$t_{\tau(i),\pi(\tau(i))}\in\mathrm{Z}(R)$ and

\noindent$\tau(i)\in\{1,\ldots,s-1,s+1,\ldots,n\}$ for all $i\in
\{1,\ldots,s-1,s+1,\ldots,n\}$, we have%
\[
t_{\tau(1),\pi(\tau(1))}\cdots t_{\tau(s-1),\pi(\tau(s-1))}t_{\tau
(s+1),\pi(\tau(s+1))}\cdots t_{\tau(n),\pi(\tau(n))}%
\]%
\[
=t_{1,\pi(1)}\cdots t_{s-1,\pi(s-1)}t_{s+1,\pi(s+1)}\cdots t_{n,\pi(n)}.
\]
The product $t_{1,\pi(1)}\cdots t_{s-1,\pi(s-1)}t_{s+1,\pi(s+1)}\cdots
t_{n,\pi(n)}$\ can be rearranged according to the cycles of $\pi$. For each
cycle $(i,\pi(i),...,\pi^{k}(i))$ of the permutation $\pi$ (of length $k+1$
say) not containing $s$ (and hence $r$) we have a factor ("sub-product")
$t_{i,\pi(i)}t_{\pi(i),\pi^{2}(i)}\cdots t_{\pi^{k}(i),\pi^{k+1}(i)}$ of the
above product and the transitivity of $T$ gives that%
\[
t_{i,\pi(i)}t_{\pi(i),\pi^{2}(i)}\cdots t_{\pi^{k}(i),\pi^{k+1}(i)}=t_{i,i}=1.
\]
The only cycle of $\pi$ containing $s$ (as well as $r$) is of the form
$(r,\pi(r),...,\pi^{l}(r)=s)$ for some $l\geq1$. The corresponding factor
("sub-product") of
\[
t_{1,\pi(1)}\cdots t_{s-1,\pi(s-1)}t_{s+1,\pi(s+1)}\cdots t_{n,\pi(n)}%
\]
does not contain $t_{\pi^{l}(r),\pi^{l+1}(r)}=t_{s,\pi(s)}=t_{s,r}$\ and the
transitivity of $T$ gives that%
\[
t_{r,\pi(r)}t_{\pi(r),\pi^{2}(r)}\cdots t_{\pi^{l-1}(r),\pi^{l}(r)}=t_{r,s}.
\]
It follows that%
\[
t_{\tau(1),\pi(\tau(1))}\cdots t_{\tau(s-1),\pi(\tau(s-1))}t_{\tau
(s+1),\pi(\tau(s+1))}\cdots t_{\tau(n),\pi(\tau(n))}=t_{r,s}%
\]
for all $\tau,\pi\in\mathrm{S}_{n}$ with $\tau(s)=s$ and $\pi(s)=r$. Thus%
\[
b(1,\tau,\pi)\cdots b(s-1,\tau,\pi)b(s+1,\tau,\pi)\cdots b(n,\tau,\pi)
\]%
\[
=t_{r,s}a_{\tau(1),\pi(\tau(1))}\cdots a_{\tau(s-1),\pi(\tau(s-1))}%
a_{\tau(s+1),\pi(\tau(s+1))}\cdots a_{\tau(n),\pi(\tau(n))},
\]
whence $\delta(a_{r,s}^{\ast})=$%
\[
t_{r,s}\!\underset{\tau,\pi}{\sum}\!\mathrm{sgn}(\pi)a_{\tau(1),\pi(\tau
(1))}\!\cdots\!a_{\tau(s-1),\pi(\tau(s-1))}\!a_{\tau(s+1),\pi(\tau
(s+1))}\!\cdots\!a_{\tau(n),\pi(\tau(n))}\!=\!t_{r,s}a_{r,s}^{\ast}%
\]
follows. $\square$

\bigskip

\noindent\textbf{4.2. Corollary.}\textit{ Let }$\delta:R\longrightarrow
R$\textit{ be an endomorphism and }$T=[t_{i,j}]$\textit{ be a transitive
matrix in }$\mathrm{M}_{n}(\mathrm{Z}(R))$\textit{. If }$A\in\mathrm{M}%
_{n}(R,\delta,T)$\textit{ is a supermatrix, then we have}%
\[
\mathrm{rdet}_{(k)}(A),\mathrm{ldet}_{(k)}(A)\in\mathrm{Fix}(\delta)
\]
\textit{for all }$k\geq1$\textit{. In particular }$\mathrm{sdet}%
(A)=\mathrm{rdet}_{(1)}(A)=\mathrm{ldet}_{(1)}(A)\in\mathrm{Fix}(\delta
)$\textit{.}

\bigskip

\noindent\textbf{Proof.} The repeated application of Theorem 4.1 gives that
the recursion $P_{1}=A^{\ast}$ and $P_{k+1}=(AP_{1}\cdots P_{k})^{\ast}$
starting from a supermatrix $A\in\mathrm{M}_{n}(R,\delta,T)$ defines a
sequence $(P_{k})_{k\geq1}$ in $\mathrm{M}_{n}(R,\delta,T)$. Since
$\mathrm{rdet}_{(k)}(A)=\mathrm{tr}(AP_{1}\cdots P_{k})$ is the sum of the
diagonal entries of the product supermatrix $AP_{1}\cdots P_{k}\in
\mathrm{M}_{n}(R,\delta,T)$ and each diagonal entry of a supermatrix (in
$\mathrm{M}_{n}(R,\delta,T)$) is in $\mathrm{Fix}(\delta)$, the proof is
complete. The poof of $\mathrm{ldet}_{(k)}(A)\in\mathrm{Fix}(\delta)$ is
similar. $\square$

\bigskip

Let $R[z]$ denote the ring of polynomials of the single commuting
indeterminate $z$, with coefficients in $R$. The $k$-th right (left)
characteristic polynomial of $A$ is the $k$-th right (left) determinant of the
$n\times n$ matrix $zI_{n}-A$ in $\mathrm{M}_{n}(R[z])$:%
\[
p_{A,k}(z)=\mathrm{rdet}_{(k)}(zI_{n}-A)\text{ and }q_{A,k}(z)=\mathrm{ldet}%
_{(k)}(zI_{n}-A).
\]
Notice that $p_{A,k}(z)$ is of the following form:%
\[
p_{A,k}(z)=\lambda_{0}^{(k)}+\lambda_{1}^{(k)}z+\cdots+\lambda_{n^{k}-1}%
^{(k)}z^{n^{k}-1}+\lambda_{n^{k}}^{(k)}z^{n^{k}},
\]
where $\lambda_{0}^{(k)},\lambda_{1}^{(k)},\ldots,\lambda_{n^{k}-1}%
^{(k)},\lambda_{n^{k}}^{(k)}\in R$ and $\lambda_{n^{k}}^{(k)}=n\left\{
(n-1)!\right\}  ^{1+n+n^{2}+\cdots+n^{k-1}}$.

\bigskip

\noindent\textbf{4.3. Corollary.}\textit{ Let }$\delta:R\longrightarrow
R$\textit{ be an endomorphism and }$T=[t_{i,j}]$\textit{ be a transitive
matrix in }$\mathrm{M}_{n}(\mathrm{Z}(R))$.\textit{ If }$A\in\mathrm{M}%
_{n}(R,\delta,T)$\textit{ is a supermatrix, then we have}%
\[
p_{A,k}(z),q_{A,k}(z)\in\mathrm{Fix}(\delta)[z]
\]
\textit{for all }$k\geq1$\textit{. In other words, the coefficients of the
right }$p_{A,k}(z)=\mathrm{rdet}_{(k)}(zI_{n}-A)$\textit{ and left }%
$q_{A,k}(z)=\mathrm{ldet}_{(k)}(zI_{n}-A)$\textit{ characteristic polynomials
are in }$\mathrm{Fix}(\delta)$\textit{.}

\bigskip

\noindent\textbf{Proof.} Now $\delta:R\longrightarrow R$ can be extended to an
endomorphism $\delta_{z}:R[z]\longrightarrow R[z]$\ of the polynomial ring
(algebra): for $r_{0},r_{1}\ldots,r_{m}\in R$ take%
\[
\delta_{z}(r_{0}+r_{1}z+\cdots+r_{m}z^{m})=\delta(r_{0})+\delta(r_{1}%
)z+\cdots+\delta(r_{m})z^{m}.
\]
Since $T$ can be considered as a transitive matrix over $\mathrm{Z}%
(R[z])=\mathrm{Z}(R)[z]$ and $zI_{n}-A\in\mathrm{M}_{n}(R[z],\delta_{z},T)$,
Corollary 4.2 gives that $p_{A,k}(z)=\mathrm{rdet}_{(k)}(zI_{n}-A)$ and
$q_{A,k}(z)=\mathrm{ldet}_{(k)}(zI_{n}-A)$ are in $\mathrm{Fix}(\delta
_{z})=\mathrm{Fix}(\delta)[z]$. $\square$

\bigskip

\noindent\textbf{4.4. Theorem.}\textit{ Let }$\delta:R\longrightarrow
R$\textit{ be an endomorphism and }$T=[t_{i,j}]$\textit{ be a transitive
matrix in }$\mathrm{M}_{n}(\mathrm{Z}(R))$. \textit{If }$R$\textit{ satisfies
the polynomial identity}%
\[
\lbrack\lbrack\lbrack\ldots\lbrack\lbrack x_{1},x_{2}],x_{3}],\ldots
],x_{k}],x_{k+1}]=0
\]
\textit{(}$R$\textit{ is Lie nilpotent of index }$k$\textit{) and }%
$A\in\mathrm{M}_{n}(R,\delta,T)$\textit{ is a supermatrix, then a right
Cayley-Hamilton identity}%
\[
(A)p_{A,k}=I_{n}\lambda_{0}^{(k)}+A\lambda_{1}^{(k)}+\cdots+A^{n^{k}-1}%
\lambda_{n^{k}-1}^{(k)}+A^{n^{k}}\lambda_{n^{k}}^{(k)}=0
\]
\textit{holds, where the coefficients }$\lambda_{i}^{(k)}$\textit{, }$0\leq
i\leq n^{k}$\textit{\ of }$p_{A,k}(z)=\mathrm{rdet}_{(k)}(zI_{n}-A)$\textit{
are in }$\mathrm{Fix}(\delta)$\textit{. If }$\lambda_{n^{k}}^{(k)}=n\left\{
(n-1)!\right\}  ^{1+n+n^{2}+\cdots+n^{k-1}}$\textit{ is invertible in }%
$R$\textit{ and }$\mathrm{Fix}(\delta)\subseteq\mathrm{Z}(R)$\textit{, then
the above identity provides the integrality of }$\mathrm{M}_{n}(R,\delta
,T)$\textit{ over }$\mathrm{Z}(R)$\textit{ (of degree }$n^{k}$\textit{).}

\bigskip

\noindent\textbf{Proof.} Since one of the main results of [S1] is that%
\[
(A)p_{A,k}=I_{n}\lambda_{0}^{(k)}+A\lambda_{1}^{(k)}+\cdots+A^{n^{k}-1}%
\lambda_{n^{k}-1}^{(k)}+A^{n^{k}}\lambda_{n^{k}}^{(k)}=0
\]
holds for $A\in\mathrm{M}_{n}(R)$, Corollary 4.3 can be used. $\square$

\bigskip

The combination of Theorems 3.2 and 4.4 yields the following.

\bigskip

\noindent\textbf{4.5. Theorem.}\textit{ Let }$R$\textit{ be a Lie nilpotent
algebra of index }$k\geq1$\textit{ over a field }$K$\textit{ of characteristic
zero and let }$e\in K$\textit{ be a primitive }$n$\textit{-th root of unity.
If }$\delta:R\longrightarrow R$\textit{ is a }$K$\textit{-automorphism
(}$K\subseteq\mathrm{Fix}(\delta)$\textit{)\ with }$\delta^{n}=\mathrm{id}%
_{R}$\textit{, then }$R$\textit{ is right (and left) integral over
}$\mathrm{Fix}(\delta)$\textit{ of degree }$n^{k}$\textit{. In other words,
for any }$r\in R$\textit{ we have}%
\[
c_{0}^{\prime}+rc_{1}^{\prime}+\cdots+r^{n^{k}-1}c_{n^{k}-1}^{\prime}%
+r^{n^{k}}=0=c_{0}^{\prime\prime}+c_{1}^{\prime\prime}r+\cdots+c_{n^{k}%
-1}^{\prime\prime}r^{n^{k}-1}+r^{n^{k}}%
\]
\textit{for some }$c_{i}^{\prime},c_{i}^{\prime\prime}\in\mathrm{Fix}(\delta
)$\textit{, }$0\leq i\leq n^{k}-1$\textit{.}

\bigskip

\noindent\textbf{Proof}. Let $P^{(e)}=[p_{i,j}]$ be the same $n\times
n$\ transitive matrix (with $p_{i,j}=e^{i-j}$, $1\leq i,j\leq n$) as in
Proposition 3.4. The application of Theorem 3.2 provides an embedding
$\overline{\delta}:R\longrightarrow\mathrm{M}_{n}(R,\delta,P^{(e)})$ such that
$\overline{\delta}(cr)=c\overline{\delta}(r)$ and $\overline{\delta
}(rc)=\overline{\delta}(r)c$ for all $r\in R$ and $c\in\mathrm{Fix}(\delta)$.
Theorem 4.4 ensures that
\[
I_{n}\lambda_{0}^{(k)}+\overline{\delta}(r)\lambda_{1}^{(k)}+\cdots
+(\overline{\delta}(r))^{n^{k}-1}\lambda_{n^{k}-1}^{(k)}+(\overline{\delta
}(r))^{n^{k}}\lambda_{n^{k}}^{(k)}=0
\]
holds, where the coefficients of the $k$-th right characteristic polynomial%
\[
p_{\overline{\delta}(r),k}(z)=\mathrm{rdet}_{(k)}(zI_{n}-\overline{\delta
}(r))=\lambda_{0}^{(k)}+\lambda_{1}^{(k)}z+\cdots+\lambda_{n^{k}-1}%
^{(k)}z^{n^{k}-1}+\lambda_{n^{k}}^{(k)}z^{n^{k}}%
\]
of $\overline{\delta}(r)\in\mathrm{M}_{n}(R,\delta,P^{(e)})$ are in
$\mathrm{Fix}(\delta)$. Since $\lambda_{n^{k}}^{(k)}=n\left\{  (n-1)!\right\}
^{1+n+n^{2}+\cdots+n^{k-1}}$ is invertible in $K$, for $c_{i}^{\prime}%
=\lambda_{i}^{(k)}\cdot\left(  \lambda_{n^{k}}^{(k)}\right)  ^{-1}%
\in\mathrm{Fix}(\delta)$, $0\leq i\leq n^{k}-1$ we have%
\[
\overline{\delta}(c_{0}^{\prime}+rc_{1}^{\prime}+\cdots+r^{n^{k}-1}c_{n^{k}%
-1}^{\prime}+r^{n^{k}})\!=\!I_{n}c_{0}^{\prime}+\overline{\delta}%
(r)c_{1}^{\prime}+\cdots+(\overline{\delta}(r))^{n^{k}-1}c_{n^{k}-1}^{\prime
}+(\overline{\delta}(r))^{n^{k}}\!=\!0.
\]
Thus $\ker(\overline{\delta})=\{0\}$ gives the desired right integrality (the
case of left integrality is similar). $\square$

\newpage

\noindent5. SUPERMATRIX\ ALGEBRAS OVER\ THE\ GRASSMANN\ ALGEBRA

\bigskip

A $\mathbb{Z}_{2}$-grading of a ring $R$ is a pair $(R_{0},R_{1})$, where
$R_{0}$ and $R_{1}$ are additive subgroups of $R$ such that $R=R_{0}\oplus
R_{1}$ is a direct sum and $R_{i}R_{j}\subseteq R_{i+j}$ for all
$i,j\in\{0,1\}$ and $i+j$ is taken modulo $2$. The relation $R_{0}%
R_{0}\subseteq R_{0}$ ensures that $R_{0}$ is a subring of $R$. Now any
element $r\in R$ can be uniquely written as $r=r_{0}+r_{1}$, where $r_{0}\in
R_{0}$ and $r_{1}\in R_{1}$. It is easy to see that the existence of $1\in R$
implies that $1\in R_{0}$. The function $\rho:R\longrightarrow R$ defined by
$\rho(r_{0}+r_{1})=r_{0}-r_{1}$ is an automorphism of $R$\ with $\rho
^{2}=\mathrm{id}_{R}$.

\bigskip

\noindent\textbf{5.1. Example.} A certain supermatrix algebra $\mathrm{M}%
_{n,d}(R)$ is considered in [S2]. In view of $\mathrm{Fix}(\rho)=R_{0}$ and
$\rho(r_{0}+r_{1})=-(r_{0}+r_{1})\Longleftrightarrow r_{0}=0$, it is
straightforward to see that $\mathrm{M}_{n,d}(R)=\mathrm{M}_{n}(R,\rho
,P(d,n))$, where $P(d,n)$ is the $n\times n$ blow up matrix in 2.4 with $m=n$.

The Grassmann (exterior) algebra%
\[
E=K\left\langle v_{1},v_{2},\ldots,v_{i},\ldots\mid v_{i}v_{j}+v_{j}%
v_{i}=0\text{ for all }1\leq i\leq j\right\rangle
\]
over a field $K$\ (of characteristic zero) generated by the infinite sequence
of anticommutative indeterminates $(v_{i})_{i\geq1}$ is a typical example of a
$\mathbb{Z}_{2}$-graded algebra. Using the well known $\mathbb{Z}_{2}$-grading
$E=E_{0}\oplus E_{1}$ and the corresponding automorphism $\varepsilon
:E\longrightarrow E$, we obtain the classical supermatrix algebra
$\mathrm{M}_{n,d}(E)=\mathrm{M}_{n}(E,\varepsilon,P(d,n))$.

Since $\varepsilon^{2}=\mathrm{id}_{E}$ implies that $(\varepsilon_{n}%
)^{2}=\mathrm{id}_{\mathrm{M}_{n}(E)}$, we can take $R=\mathrm{M}_{n}(E)$ and
$\delta=\varepsilon_{n}$ in Corollary 3.5. Thus we obtain the well known
embedding%
\[
\overline{(\varepsilon_{n})}:\mathrm{M}_{n}(E)\longrightarrow\mathrm{M}%
_{2}(\mathrm{M}_{n}(E),\varepsilon_{n},P)\cong\mathrm{M}_{2n}(E,\varepsilon
,P(n,2n))=\mathrm{M}_{2n,n}(E)
\]
of $\mathrm{M}_{n}(E)$\ into the supermatrix algebra $\mathrm{M}_{2n,n}(E)$.

In view of $\mathrm{Fix}(\varepsilon)=E_{0}=\mathrm{Z}(E)$, the application of
Theorem 4.4 gives that $\mathrm{M}_{n,d}(E)$ is integral over $E_{0}$ of
degree $n^{2}$ (see Theorem 3.3 in [S2]). Thus the main results of [S2]
directly follow from 4.1, 4.2, 4.3 and 4.4. $\square$

\bigskip

We note that the T-ideal of the polynomial identities (with coefficients in
$K$) satisfied by $\mathrm{M}_{n,d}(E)$ plays an important role in Kemer's
classification of the T-prime T-ideals (see [K]).

\bigskip

\noindent\textbf{5.2. Example.} For an integer $k\geq0$ let%
\[
E(k)=\underset{1\leq i_{1}<\ldots<i_{k}}{{\LARGE \oplus}}Kv_{i_{1}}\cdots
v_{i_{k}}%
\]
denote the $k$-homogeneous component of $E$, i.e. the $K$-linear span of all
products $v_{i_{1}}\cdots v_{i_{k}}$ of length $k$ ($E(0)=K$). If $e\in K$ is
a primitive $n$-th root of unity, then define an automorphism $\rho
_{e}:E\longrightarrow E$ as follows%
\[
\rho_{e}(h(v_{1},v_{2},\ldots,v_{k},\ldots))=h(ev_{1},ev_{2},\ldots
,ev_{k},\ldots),
\]
where each generator $v_{k}$ in $h\in E$ is replaced by $ev_{k}$.

Consider the supermatrix algebra $\mathrm{M}_{n}(E,\rho_{e},P^{(e)})$, where
$P^{(e)}=[p_{i,j}]$ is the same $n\times n$ transitive matrix (with
$p_{i,j}=e^{i-j}$, $1\leq i,j\leq n$) as in Proposition 3.4. If $A=[h_{i,j}]$
is a supermatrix in $\mathrm{M}_{n}(E,\rho_{e},P^{(e)})$, then we have
$\rho_{e}(h_{i,j})=e^{i-j}h_{i,j}$ for all $1\leq i,j\leq n$. Since for an
integer $0\leq m\leq n-1$%
\[
(ev_{i_{1}})\cdots(ev_{i_{k}})=e^{k}v_{i_{1}}\cdots v_{i_{k}}=e^{m}v_{i_{1}%
}\cdots v_{i_{k}}%
\]
holds if and only if $k-m$ is divisible by $n$, we obtain that%
\[
E_{m,n}=\{h\in E\mid\rho_{e}(h)=e^{m}h)\}=\underset{u=0}{\overset{\infty
}{{\LARGE \oplus}}}E(m+nu)
\]
and%
\[
E_{-m,n}=\{h\in E\mid\rho_{e}(h)=e^{-m}h)\}=E_{n-m,n}.
\]
Thus the shape of $\mathrm{M}_{n}(E,\rho_{e},P^{(e)})$ is the following:%
\[
\mathrm{M}_{n}(E,\rho_{e},P^{(e)})=\left[
\begin{array}
[c]{ccccc}%
E_{0,n} & E_{-1,n} & \cdots & E_{-(n-2),n} & E_{-(n-1),n}\\
E_{1,n} & E_{0,n} & E_{-1,n} & \ddots & E_{-(n-2),n}\\
\vdots & E_{1,n} & \ddots & \ddots & \vdots\\
E_{n-2,n} & \ddots & \ddots & E_{0,n} & E_{-1,n}\\
E_{n-1,n} & E_{n-2,n} & \cdots & E_{1,n} & E_{0,n}%
\end{array}
\right]  .
\]
Since $(\rho_{e})^{n}=\mathrm{id}_{E}$, Theorem 3.2 provides an embedding
$\overline{\rho_{e}}:E\longrightarrow\mathrm{M}_{n}(E,\rho_{e},P^{(e)})$. In
view of the Lie nilpotency of $E$ (of index $2$), the application of Theorem
4.4 gives that any matrix $A\in\mathrm{M}_{n}(E,\rho_{e},P^{(e)})$ satisfies a
right Cayley-Hamilton identity of degree $n^{2}$ with (right) coefficients
from $\mathrm{Fix}(\rho_{e})=E_{0,n}$. If $n$ is even, then $E_{0,n}\subseteq
E_{0}=\mathrm{Z}(E)$ is a central subalgebra and the coefficients in the
mentioned (right-left) Cayley-Hamilton identity are central. The application
of Theorem 4.5 gives that $E$ is right (and left) integral over $E_{0,n}$ of
degree $n^{2}$.

If $n=2$ and $e=-1$, then $E_{0,2}=E_{0}$ is the even and $E_{1,2}=E_{1}$ is
the odd part of the Grassmann algebra $E$ and $\mathrm{M}_{2}(E,\rho
_{-1},P^{(-1)})=\mathrm{M}_{2}(E,\varepsilon,P)$ (for $P$ see 2.4). $\square$

\bigskip

\noindent\textbf{5.3. Example.} For $g\in E$, let $\sigma:E\longrightarrow E$
be the following map:%
\[
\sigma(g)=(1+v_{1})g(1-v_{1}).
\]
Clearly, $1-v_{1}=(1+v_{1})^{-1}$ implies that $\sigma$ is a conjugate
automorphism of $E$. The blow up $Q(d,n)$\ of the transitive matrix%
\[
Q=\left[
\begin{array}
[c]{cc}%
1 & 1+v_{1}v_{2}\\
1-v_{1}v_{2} & 1
\end{array}
\right]
\]
is in $\mathrm{M}_{n}(E_{0})$ (notice that $E_{0}=\mathrm{Z}(E)$). Thus we can
form the supermatrix algebra $\mathrm{M}_{n}(E,\sigma,Q(d,n))$. The block
structure of a supermatrix $A\in\mathrm{M}_{n}(E,\sigma,Q(d,n))$ is the
following%
\[
A=\left[
\begin{array}
[c]{cc}%
A_{1,1} & A_{1,2}\\
A_{2,1} & A_{2,2}%
\end{array}
\right]  ,
\]
where the square blocks $A_{1,1}$ and $A_{2,2}$ are of sizes $d\times d$ and
$(m-d)\times(m-d)$ and the rectangular blocks $A_{1,2}$ and $A_{2,1}$ are of
sizes $d\times(m-d)$ and $(m-d)\times d$. The entries of $A_{1,1}$ and
$A_{2,2}$ are in%
\[
\mathrm{Fix}(\sigma)\!=\!\{g\in E\!\mid\!(1\!+\!v_{1})g(1\!-\!v_{1}%
)\!=\!g\}\!=\!\{g\in E\!\mid\!v_{1}g\!-\!gv_{1}\!=\!0\}\!=\!\mathrm{Cen}%
(v_{1})\!=\!E_{0}\!+\!E_{0}v_{1},
\]
where $\mathrm{Cen}(v_{1})$ denotes the centralizer of $v_{1}$. The entries of
$A_{1,2}$ are in%
\[
\Omega_{1,2}=\{g\in E\mid(1+v_{1})g(1-v_{1})=(1+v_{1}v_{2})g\}=\{g\in E\mid
v_{1}g-gv_{1}=v_{1}v_{2}g\}
\]%
\[
=\{g_{0}+g_{1}\mid g_{0}\in E_{0},g_{1}\in E_{1},v_{1}g_{1}-g_{1}v_{1}%
=v_{1}v_{2}g_{0}\text{ and }v_{1}v_{2}g_{1}=0\}
\]%
\[
=\{g_{0}+g_{1}\mid g_{0}\in E_{0},2g_{1}-v_{2}g_{0}\in E_{0}v_{1}\}\subseteq
E_{0}+E_{0}v_{1}+E_{0}v_{2}.
\]
The entries of $A_{2,1}$ are in%
\[
\Omega_{2,1}=\{g\in E\mid(1+v_{1})g(1-v_{1})=(1-v_{1}v_{2})g\}=\{g\in E\mid
v_{1}g-gv_{1}=-v_{1}v_{2}g\}
\]%
\[
=\{g_{0}+g_{1}\mid g_{0}\in E_{0},g_{1}\in E_{1},v_{1}g_{1}-g_{1}v_{1}%
=-v_{1}v_{2}g_{0}\text{ and }v_{1}v_{2}g_{1}=0\}
\]%
\[
=\{g_{0}+g_{1}\mid g_{0}\in E_{0},2g_{1}+v_{2}g_{0}\in E_{0}v_{1}\}\subseteq
E_{0}+E_{0}v_{1}+E_{0}v_{2}%
\]
As a consequence, we obtain that the shape of $\mathrm{M}_{n}(E,\sigma
,Q(d,n))$ is the following:%
\[
\mathrm{M}_{n}(E,\sigma,Q(d,n))=\left[
\begin{array}
[c]{cc}%
E_{0}+E_{0}v_{1} & \Omega_{1,2}\\
\Omega_{2,1} & E_{0}+E_{0}v_{1}%
\end{array}
\right]
\]
with diagonal blocks of sizes $d\times d$ and $(m-d)\times(m-d)$. In view of
the Lie nilpotency of $E$ (of index $2$), the application of Theorem 4.4 gives
that any matrix $A\in\mathrm{M}_{n}(E,\sigma,Q(d,n))$ satisfies a right
Cayley-Hamilton identity of degree $n^{2}$ with (right) coefficients from
$E_{0}+E_{0}v_{1}$. $\square$

\bigskip

\noindent REFERENCES

\bigskip

\noindent\lbrack Do] M. Domokos, \textit{Cayley-Hamilton theorem for }%
$2\times2$\textit{\ matrices over the Grassmann algebra}, J. Pure Appl.
Algebra 133 (1998), 69-81.

\noindent\lbrack K] A. R. Kemer,\textit{\ Ideals of Identities of Associative
Algebras,} Translations of Math. Monographs, Vol. 87 (1991), AMS, Providence,
Rhode Island.

\noindent\lbrack S1] J. Szigeti, \textit{New determinants and the
Cayley-Hamilton theorem for matrices over Lie nilpotent rings}, Proc. Amer.
Math. Soc. 125 (1997), 2245-2254.

\noindent\lbrack S2] J. Szigeti, \textit{On the characteristic polynomial of
supermatrices}, Israel Journal of Mathematics Vol. 107 (1998), 229-235.

\noindent\lbrack S3] J. Szigeti, \textit{Cayley-Hamilton theorem for matrices
over an arbitrary ring}, Serdica Math. J. 32 (2006), 269-276.

\noindent\lbrack SvW] J. Szigeti and L. van Wyk:\textit{ Determinants for
}$n\times n$\textit{ matrices and the symmetric Newton formula in the
}$3\times3$\textit{ case,} Linear and Multilinear Algebra, Vol. 62, No. 8
(2014), 1076-1090.

\end{document}